\documentclass[preprint,11pt]{elsarticle}

\newtheorem{theorem}{Theorem}[section]
\newtheorem{corollary}{Corollary}[section]
\newtheorem{lemma}{Lemma}[section]

\usepackage{amssymb}

\begin{document}

\begin{frontmatter}

\title{On the Wassertein distance for a martingale central limit theorem}
\author{Xiequan Fan$^*$,  Xiaohui Ma}
 \cortext[cor1]{\noindent Corresponding author. \\
  \mbox{ \ \ \ } \textit{E-mail}: fanxiequan@hotmail.com (X. Fan). }
\address{Center for Applied Mathematics, Tianjin University, 300072 Tianjin,  China}

\begin{abstract}
We prove an upper bound on the Wassertein distance between normalized martingales and the standard normal random variable, which extends a result  of R\"ollin \cite{AR18}. The proof is based on a method of Bolthausen \cite{EB82}.
\end{abstract}

\begin{keyword} Martingales; Central limit theorem;  Wassertein metric
\vspace{0.3cm}
\MSC   60G42; 60E15; 60F25
\end{keyword}

\end{frontmatter}


\section{Introduction and main result}
Let $n\geq1$.
Assume that $\mathbf{X}=\left(X_{i}\right)_{1\leq i\leq n}$ is a sequence of martingale differences defined on probability space $(\Omega, \mathcal{F}, \mathbf{P})$, $\mathcal{F}_0=\{\emptyset, \Omega \} $ and that $X_{i}$ is adapted to the $\sigma$-fields $\mathcal{F}_{i}, i=1,2, \ldots, n$,  where $\mathcal{F}_{i}$ is the $\sigma$-algebra generated by $X_{1}, X_{2}, \ldots, X_{i}$.
In other words, $X_{i}$ satisfies $\mathbf{E}\left[X_{i} | \mathcal{F}_{i-1}\right]=0$. Let $M_{n}$ denote the class of all such sequences of length $n$.
If $\mathbf{X} \in M_{n}$, write
$$S_{n}=X_{1}+X_{2}+\ldots+X_{n},$$
$$\sigma_{i}^{2}=\mathbf{E}\left[X_{i}^{2} |\mathcal{F}_{i-1}\right], i=1,2, \ldots, n,$$
$$s_{n}^{2}=\sum_{i=1}^{n} \mathbf{E}X_{i}^{2},$$
$$V_{n}^{2}=\frac{\sum_{i=1}^{n} \sigma_{i}^{2}}{s_{n}^{2}},$$
$$\|\mathbf{X}\|_{p}=\max _{1 \leq i \leq n}\left\|X_{i}\right\|_{p},\ \ p \in[1,+\infty],$$
where $\left\|X_{i}\right\|_{p}= \left( \mathbf{E}\left|X_{i}\right|^{p}\right)^{1 / p}$ for $1\leq p < \infty$ and
$\left\|X_{i}\right\|_{\infty}=\inf\left\{a:\mathbf{P}(|X_{i}|\leq a)=1\right\}$.

According to the martingale central limit theory, it is well known that if, as $ n \rightarrow\infty,$
$$
V_{n}^{2}\stackrel{\mathbf{P}}{\longrightarrow} 1
$$
and the ¡°conditional Lindeberg condition¡±
 $$
\sum_{i=1}^{n} \mathbf{E}\left[X_{i}^{2} \mathbf{1}_{\left\{\left|X_{i}\right| \geq \varepsilon\right\}} | \mathcal{F}_{i-1}\right] \stackrel{\mathbf{P}}{\longrightarrow} 0 \quad \textrm {for each }\ \varepsilon>0
$$
are satisfied, then  ${S_{n}/s_{n}}$ converges to the standard Gaussian random variable in distribution, that is,
$$
\forall t \in \mathbf{R}, \quad \mathbf{P}\left(S_{n} / s_{n} \leq t\right) {\longrightarrow} \Phi(t), \quad \textrm { as }\  n \rightarrow\infty,
$$
where $\Phi(t)=(2 \pi)^{-1 / 2} \int_{-\infty}^{t} \mathrm{e}^{-x^{2} / 2} \mathrm{d} x$.
Define
$$
D\left(S_{n} / s_{n}\right)=\sup _{t}|\mathbf{P}(S_{n} / s_{n} \leq t)-\Phi(t)|,
$$
then $$ D\left(S_{n} / s_{n}\right){\longrightarrow} 0, \quad \textrm { as }\  n \rightarrow\infty.
$$

The Kolmogorov distance in central limit theorem for martingales has been intensely studied under various conditions.
For instance, we recall the following result due to Heyde and Brown \cite{Heyde70}: for a constant $ p \in(1,2] $, they proved that
\begin{eqnarray}\label{Heyde}
D\left(S_{n} / s_{n}\right) \leq C_{p}\bigg(s_{n}^{-4p}\left\|V_{n}^{2}-1\right\|_{p}^{p}+s_{n}^{-2p}(\sum_{i=1}^{n} \mathbf{E}\left|X_{i}\right|^{2p})\bigg)^{1 /(2p+1)},
\end{eqnarray}
where, here and after, $C_{p}$ is a constant depending only on $p$.
Bolthausen \cite{EB82} proved that if $\|\mathbf{X}\|_{\infty} \leq \gamma$ a.s., then
\begin{eqnarray}\label{Bolth}
D\left(S_{n} / s_{n}\right) \leq C_{\gamma}\bigg(\frac{n \log n}{s_{n}^{3}}+\min \left\{\left\|V_{n}^{2}-1\right\|_{\infty}^{1 / 2},\left\|V_{n}^{2}-1\right\|_{1}^{1 / 3}\right\}\bigg).
 \end{eqnarray}
 Using a modification of the method developed by Bolthausen, Haeusler \cite{Hae88} gave an extension of (\ref{Heyde}) to all $p>1$; El Machkouri and Ouchti \cite{Mac07} replaced the term ${n \log n}/s_{n}^{3}$ of (\ref{Bolth}) by ${\max _{1 \leq i \leq n}\gamma_{i}\log n}/{s_{n}}$ to a large class of  martingale difference sequences satisfying
\begin{eqnarray}
\mathbf{E}\left[|X_{i}|^{3} | \mathcal{F}_{i-1}\right] \leq \gamma_{i}\mathbf{E}\left[|X_{i}|^{2} | \mathcal{F}_{i-1}\right] a.s.,
\end{eqnarray}
where $\gamma_{i}$, $i=1,2, \ldots, n $ are constants. Following Bolthausen again, Mourrat \cite{JCM13} has extended the term $\min \left\{\left\|V_{n}^{2}-1\right\|_{\infty}^{1 / 2},\left\|V_{n}^{2}-1\right\|_{1}^{1 / 3}\right\}$ of (\ref{Bolth}) to the more general term $\left(\left\|V_{n}^{2}-1\right\|_{p}^{p}+s_{n}^{-2 p}\right)^{1 /(2p+1)}$, for $p \geq 1.$
Recently, with the methods of Grama and Hausler \cite{Gra00} and Bolthausen \cite{EB82}  (see also Fan et al.\ \cite{FGLS19}), Fan \cite {F19} proved that if there exist two positive numbers $\rho$ and $\gamma$, such that
\begin{eqnarray}
\mathbf{E}\left[|X_{i}|^{2+\rho} | \mathcal{F}_{i-1}\right] \leq \gamma^{\rho}\mathbf{E}\left[|X_{i}|^{2} | \mathcal{F}_{i-1}\right] a.s.
\end{eqnarray}
for all $i=1,2, \ldots, n$, then for $p \geq 1,$
\begin{eqnarray}
D\left(S_{n} / s_{n}\right) \leq C_{p,\gamma,\rho}\bigg(\alpha_{n}+\Big(\left\|V_{n}^{2}-1\right\|_{p}^{p}+\frac{1}{s_{n}^{2p}}\mathbf{E} \max _{1 \leq i \leq n}\left|X_{i}\right|^{2p}\Big)^{1/(2p+1)}\bigg),
\end{eqnarray}
where
$$\alpha_{n}  = \left\{ \begin{array}{l}
\frac{1}{s_{n}^{\rho}} \ \ \ \ \ \  \ \ \ \ {\textrm{ if }} \rho \in (0,1) ,\\ \\
\frac{1}{s_{n}}\log s_{n} \ \ \  {\textrm{if }} \rho \geq 1.
\end{array} \right.$$

Despite the Kolmogorov distance in central limit theorem has been intensely studied, research of bounds with respect to the Wassertein distance are rare.
To the best of our knowledge, we only aware the articles of  Van Dung et al. \cite{LVD14} and R\"ollin \cite{AR18}.
Denote by $d_{w}\big(S_{n} / s_{n}\big)=\int_{-\infty}^{+\infty}\big|\mathbf{P}\big(S_{n} / s_{n}\leq x\big)-\Phi\big(x\big)\big|dx$
the Wassertein distance between the distributions of $S_{n}/s_{n}$ and the standard normal random variable.
Van Dung et al. \cite{LVD14} extended Mourrat's bounds to the $L^{1}$-bound in the mean central limit theorem, that is, if $\|\mathbf{X}\|_{\infty} \leq \gamma$ a.s., then for $p > 1/2$,
\begin{eqnarray}\label{v}
d_{w}\big(S_{n} / s_{n}\big) \leq  C_{p, \gamma}\bigg(\frac{ n \log n}{s_{n}^{3}}+\left(\left\|V_{n}^{2}-1\right\|_{p}^{p}+s_{n}^{-2 p}\right)^{1 / 2 p}\bigg).
\end{eqnarray}
 R\"ollin \cite{AR18} provided a new proof of already known result by using combination of both Lindeberg's and Stein's methods. His result states that if $V_{n}^{2}=1$ a.s., then for any $a\geq0$,
 \begin{eqnarray}\label{r}
 d_{w}\big(S_{n} / s_{n}\big) \leq \frac{3}{s_{n}}\sum_{i=1}^{n} \mathbf{E}\frac{\left|X_{i}\right|^{3}}{{\rho}_{i}^{2}+a^2}+\frac{2a}{s_{n}},
 \end{eqnarray}
where $\rho_{i}^{2}=\sum_{k=i}^{n}\sigma_{k}^{2}$.
The aim of this article is to extend  (\ref{r})  by relaxing the condition $V_{n}^{2}=1$ a.s. to $||V_{n}^{2}-1||_{p}{\longrightarrow} 0$ for any $p\geq 1$.

The following theorem is our main result.
\begin{theorem}\label{theorem}
For any $p\geq 1$, there exists a constant $C_{p}>0$ such that
\begin{eqnarray}\label{thm inequa}
d_{w}\big(S_{n} / s_{n}\big)
\leq C_{p}\bigg(\frac{1}{s_{n}}\Big(\sum_{i=1}^{n} \mathbf{E}\left|X_{i}\right|^{3}\Big)^{1/3}+\Big(\left\|V_{n}^{2}-1\right\|_{p}^{p}+\frac{1}{s_{n}^{2p}}\mathbf{E} \max _{1 \leq i \leq n}\left|X_{i}\right|^{2p}\Big)^{1/2p}\bigg).
\end{eqnarray}
\end{theorem}

Notice that if $\|\mathbf{X}\|_{\infty} \leq \gamma$ a.s.,
then the order of the term $\Big(\|V_{n}^{2}-1\|_{p}^{p}+\frac{1}{s_{n}^{2p}}\mathbf{E} \max _{1 \leq i \leq n}|X_{i}|^{2p}\Big)^{1/2p}$ is less than
that of $\Big(\|V_{n}^{2}-1\|_{p}^{p}+s_{n}^{-2 p}\Big)^{1 / 2 p}$. Thus the second term of bound (\ref{thm inequa}) implies the second term of bound  (\ref{v}).

Let $p=3/2.$ The following corollary is an immediate consequence of Theorem \ref{theorem}.
\begin{corollary}\label{cor2}
There exists a constant $C>0$ such that
\begin{eqnarray}\label{cor2in}
d_{w}\big(S_{n} / s_{n}\big)
\leq C\bigg(\frac{1}{s_{n}}\Big(\sum_{i=1}^{n} \mathbf{E}\left|X_{i}\right|^{3}\Big)^{1/3}+\left\|V_{n}^{2}-1\right\|_{3/2}^{1/2}\bigg).
\end{eqnarray}
\end{corollary}

Though our bound may not be smaller than the bound of R\"ollin (cf.\ (\ref{r})), the advantage of (\ref{cor2in}) is that the term $\rho_{i}^{2}$ appearing in (\ref{r}) does not appear any more. Moreover, inequality (\ref{cor2in}) gives a rate of convergence under weaker condition, that is the condition $V_{n}^{2}=1$ a.s.\ has been relaxed to $\left\|V_{n}^{2}-1\right\|_{3/2}^{1/2}{\longrightarrow} 0$.

\section{Proof of Theorem \ref{theorem}}
In the proof of the Theorem \ref{theorem}, we shall use the following  lemma which is a consequence of   R\"ollin's inequality (\ref{r}) by dropping  $\rho_{i}^{2}$.
\begin{lemma}\label{lemma2}
Assume that $V_{n}^{2}=1$ a.s., there exists a constant $C>0$ such that
\begin{eqnarray}\label{rr}
d_{w}\big(S_{n} / s_{n}\big) \leq \frac{C}{s_{n}}\left( \sum_{i=1}^{n} \mathbf{E}\left|X_{i}\right|^{3}\right)^{1 / 3}.
\end{eqnarray}
\end{lemma}
\noindent\emph{Proof of Lemma \ref{lemma2}.}
From (\ref{r}), by the fact $\rho_{i}^{2} \geq 0,$ it is easy to see that
\begin{eqnarray}
 d_{w}\big(S_{n} / s_{n}\big) \leq \frac{3}{s_{n}}\sum_{i=1}^{n}\frac{ \mathbf{E}\left|X_{i}\right|^{3}}{a^2}+\frac{2a}{s_{n}}.
 \end{eqnarray}
Let $f(a)=\sum_{i=1}^{n}\frac{ \mathbf{E}\left|X_{i}\right|^{3}}{a^2}+2a.$
Clearly if $a=\left( \sum_{i=1}^{n} \mathbf{E}\left|X_{i}\right|^{3}\right)^{1 / 3},$ then $f(a)$ achieves the minimum, which gives the desired inequality (\ref{rr}).
\hfill\qed \\

Now, we are in position to prove Theorem \ref{theorem}.
The idea is to construct a new sequence of martingale differences  $\big(\hat{X_{i}},  \hat{\mathcal{F}}_{i}\big)_{1\leq i\leq N}$ based on $\mathbf{X}$  such that $\sum_{i=1}^{N} \mathbf{E}\big[\hat{X}_{i}^{2} | \hat{\mathcal{F}}_{i-1}\big]/s_{n}^{2}   =1$ a.s., where the exact value of $N$ is given later, and then apply Lemma \ref{lemma2} to $\hat{\mathbf{X}} $. Consider  the stopping time
$$
\tau=\sup \bigg\{1\leq k \leq n : \sum_{i=1}^{k}\mathbf{E} \left[X_{i}^{2} | \mathcal{F}_{i-1}\right] \leq s_{n}^{2}\bigg\}.
$$
Assume that $\varepsilon>0$. 
Let $r=\left\lfloor\frac{s_{n}^{2}-\sum_{i=1}^{\tau} \mathbf{E}\left[X_{i}^{2} | \mathcal{F}_{i-1}\right]}{\varepsilon^{2}}\right\rfloor,$ where $\left\lfloor x \right\rfloor$ stands for the largest integer not exceeding $x.$
Clearly $r \leq \left\lfloor s_{n}^{2}/\varepsilon^{2}\right\rfloor$.
Let $N=n+\left\lfloor s_{n}^{2}/\varepsilon^{2}\right\rfloor+1.$
Conditionally on $\mathcal{F}_{\tau}$, and for $\tau+1 \leq i \leq \tau+r$,   let $Y_{i}$ be independent random variables such that $\mathbf{P}(Y_{i}=\pm \varepsilon)=1 / 2.$
When $i=\tau+r+1$,  let $Y_{\tau+r+1}$ be such that
$$
\mathbf{P}\bigg( Y_{\tau+r+1}=\pm\Big(s_{n}^{2}-\sum_{i=1}^{\tau} \mathbf{E}\left[X_{i}^{2} | \mathcal{F}_{i-1}\right]-r \varepsilon^{2}\Big)^{1 / 2}\bigg)=\frac{1}{2},
$$
with the sign determined independent of everything else. Finally, if $\tau+r+1<i\leq N$, then we let $Y_{i}=0$. By the definition of $\hat{X}_i$,
we have
$$\hat{X}_i= X_i \mathbf{1}_{\{i \leq \tau  \}} + Y_i \mathbf{1}_{\{  \tau  < i \leq \tau +r \}} + Y_{\tau +r +1} \mathbf{1}_{\{ i=   \tau +r +1 \}},\ \ 1\leq i \leq N. $$
Define $\hat{\mathcal{F}}_{i}= \mathcal{F}_{i} $ for $i \leq \tau,$ $\hat{\mathcal{F}}_{i}=\sigma\{\mathcal{F}_{\tau},\  Y_{j}, \tau+1 \leq j \leq i \}$ for $\tau+1\leq i \leq N.$
%
Since that $\mathbf{1}_{\{i \leq \tau  \}}=1- \mathbf{1}_{\{  \tau \leq i-1  \}} $ is $\mathcal{F}_{i-1}$ measurable, we deduce that
\begin{eqnarray*}
\mathbf{E}[  \hat{X}_i | \hat{\mathcal{F}}_{i-1}] &=& \mathbf{E}[ X_i \mathbf{1}_{\{i \leq \tau  \}} | \hat{\mathcal{F}}_{i-1}] + \mathbf{E}[ Y_i \mathbf{1}_{\{  \tau  < i \leq \tau +r \}} | \hat{\mathcal{F}}_{i-1}]+ \mathbf{E}[  Y_{ \tau +r +1} \mathbf{1}_{\{ i=   \tau +r +1 \}}| \hat{\mathcal{F}}_{i-1}] \nonumber \\
 &=&\mathbf{E}[ X_i\mathbf{1}_{\{i \leq \tau  \}} |  \mathcal{F}_{i-1}] + \Big( \frac{\varepsilon}2 \mathbf{E}[  \mathbf{1}_{\{  \tau  < i \leq \tau +r \}} | \hat{\mathcal{F}}_{i-1}] -\frac{\varepsilon}2\mathbf{E}[   \mathbf{1}_{\{  \tau  < i \leq \tau +r \}} | \hat{\mathcal{F}}_{i-1}]  \Big) \nonumber \\
  &&+  \bigg( \frac12 \Big(s_{n}^{2}-\sum_{i=1}^{\tau} \mathbf{E}\left[X_{i}^{2} | \mathcal{F}_{i-1}\right]-r \varepsilon^{2}\Big)^{1 / 2} - \frac12\Big(s_{n}^{2}-\sum_{i=1}^{\tau} \mathbf{E}\left[X_{i}^{2} | \mathcal{F}_{i-1}\right]-r \varepsilon^{2}\Big)^{1 / 2} \bigg) \nonumber \\
  &=& 0 \ \ \ \ \textrm{a.s.}
 \end{eqnarray*}
 Thus $\hat{\mathbf{X}}=(\hat{X}_i, \hat{\mathcal{F}}_i)_{1\leq i \leq N}$ is also a martingale difference sequence.
Moreover, it holds
$$ \sum_{i=\tau+1}^{N} \mathbf{E}\big[\hat{X}_{i}^{2} | \hat{\mathcal{F}}_{i-1}\big]= \sum_{i=\tau+1}^{\tau+r} \varepsilon^{2}+ s_{n}^{2}-\sum_{i=1}^{\tau} \mathbf{E}\left[X_{i}^{2} | \mathcal{F}_{i-1}\right]-r \varepsilon^{2}=s_{n}^{2}-\sum_{i=1}^{\tau} \mathbf{E}\left[X_{i}^{2} | \mathcal{F}_{i-1}\right], $$
which implies that
$$
\sum_{i=1}^{N} \mathbf{E}\big[\hat{X}_{i}^{2} | \hat{\mathcal{F}}_{i-1}\big]=s_{n}^{2}\quad \quad a.s.
$$
Consequently, $\sum_{i=1}^{N} \mathbf{E}\big[\hat{X}_{i}^{2} | \hat{\mathcal{F}}_{i-1}\big]/s_{n}^{2}=1$ a.s., by Lemma \ref{lemma2} it is easy to see that:
\begin{eqnarray}
d_{w}\big(\hat S_{N} / s_{n}\big)
\leq \frac{C}{s_{n}}\bigg( \sum_{i=1}^{N} \mathbf{E} |\hat{X}_{i} |^{3}\bigg)^{1 / 3}.
\end{eqnarray}
According to the construction of $\hat{\mathbf{X}}$, we can easily give an estimation of $ \sum_{i=1}^{N} \mathbf{E}|\hat{X}_{i}|^{3} $ as following:
$$\sum_{i=1}^{N} \mathbf{E}|\hat{X}_{i}|^{3}=\sum_{i=1}^{\tau} \mathbf{E}\left|{X}_{i}\right|^{3}+\sum_{i=\tau +1}^{N} \mathbf{E}|\hat{X}_{i}|^{3}\nonumber \\
\leq  \sum_{i=1}^{n} \mathbf{E}\left|{X}_{i}\right|^{3}+\left(1+s_{n}^{2}/\varepsilon^{2}\right)\varepsilon^{3}.$$
Thus we get,
\begin{eqnarray}\label{1}
d_{w}\big(\hat S_{N} / s_{n}\big)  \leq \frac{C}{s_{n}}\left( \sum_{i=1}^{n} \mathbf{E}\left|X_{i}\right|^{3}+\varepsilon^{3}+s_{n}^{2}\varepsilon \right)^{1 / 3}.
\end{eqnarray}
For any $x>0,$ we have
\begin{eqnarray}\label{a}
d_{w}\big(S_{n} / s_{n}\big)  &\leq& d_{w}\big(\hat S_{N} / s_{n}\big)+\int_{-\infty}^{+\infty}\Big|\mathbf{P}\big(\hat{S}_{N} / s_{n}\leq t+x\big)-\mathbf{P}\big({S}_{n} / s_{n}\leq t\big)\Big|dt+\int_{-\infty}^{+\infty}\Big|\Phi(t+x)-\Phi(t)\Big|dt \nonumber \\
&\leq&
d_{w}\big(\hat S_{N} / s_{n}\big)+\frac{C}{x^{2p-1}s_{n}^{2p}}\mathbf{E}\big|\hat{S}_{N}-S_{n}\big|^{2p}+2x,
\end{eqnarray}
see (5.3) of Van Dung et al.\ \cite{LVD14}.
We give an estimation for the second term in the right-hand side of (\ref{a}). First we note that
$$ S_{n}-\hat{S}_{N} =\sum_{i=\tau+1}^{N}\big(X_{i}-\hat{X}_{i}\big),$$
where we put $X_{i}=0$ for $i>n$.
 As $\tau$ is a stopping time, conditionally on $\tau$, the $(X_{i}-\hat{X}_{i} )_{ \tau+1 \leq i \leq N }$ still forms a martingale difference sequence, see Mourrat \cite{JCM13} for details.
 By Burkholder's inequality, we deduce that
\begin{eqnarray}\label{aa}
{\frac{1}{C_{p}} \mathbf{E}\left[\big|\hat{S}_{N}-S_{n}\big|^{2 p}\right]}
  \leq  \mathbf{E}\left[\left(\sum_{i=\tau+1}^{N} \mathbf{E}\left[\big(X_{i}-\hat{X}_{i}\big)^{2} \Big | \mathcal{F}_{i-1}\right]\right)^{p}\, \right]+\mathbf{E}\left[\max _{\tau+1 \leq i \leq N}\left|X_{i}-\hat{X}_{i}\right|^{2 p}\right].
\end{eqnarray}
As $\mathbf{E} [X_{i} \hat{X}_{i} | \mathcal{F}_{i-1} ]=0$ for $i\geq \tau+1$, we have
\begin{eqnarray}
 \sum_{i=\tau+1}^{N} \mathbf{E}\left[\big(X_{i}-\hat{X}_{i}\big)^{2} \Big| \mathcal{F}_{i-1}\right] &=&
 \sum_{i=\tau+1}^{N} \mathbf{E}\left[X_{i}^{2} | \mathcal{F}_{i-1}\right]+\sum_{i=\tau+1}^{N} \mathbf{E}\, \big[\hat{X}_{i}^{2} | \mathcal{F}_{i-1}\big] \nonumber\\
  &=&s_{n}^{2} V_{n}^{2}+s_{n}^{2}-2 \sum_{i=1}^{\tau +1} \mathbf{E}\left[X_{i}^{2} | \mathcal{F}_{i-1}\right]+2\, \mathbf{E}\left[X_{\tau +1}^{2} | \mathcal{F}_{\tau}\right] .
\end{eqnarray}
Notice that $\sum_{i=1}^{\tau+1} \mathbf{E}\left[X_{i}^{2} | \mathcal{F}_{i-1}\right] > s_{n}^{2}.$ Hence, we get
\begin{eqnarray}
\sum_{i=\tau+1}^{N} \mathbf{E}\left[\big(X_{i}-\hat{X}_{i}\big)^{2} \Big | \mathcal{F}_{i-1}\right] \leq s_{n}^{2}V_{n}^{2}-s_{n}^{2}+2\mathbf{E}\left[X_{\tau +1}^{2} | \mathcal{F}_{\tau}\right] .
\end{eqnarray}
Using the inequality $|a+b|^{k} \leq 2^{k-1}\left(|a|^{k}+|b|^{k}\right), k\geq1$,  we have
\begin{eqnarray}\label{aaa}
\left(\sum_{i=\tau+1}^{N} \mathbf{E}\left[\big(X_{i}-\hat{X}_{i}\big)^{2} \Big | \mathcal{F}_{i-1}\right]\right)^{p}
\leq 2^{p-1}|s_{n}^{2}V_{n}^{2}-s_{n}^{2}|^p+2^{2p-1}\left(\mathbf{E}\left[X_{\tau +1}^{2} | \mathcal{F}_{\tau}\right]\right)^p.
\end{eqnarray}
For the second term in the right-hand side of (\ref{aaa}), by Jensen's inequality, we obtain
\begin{eqnarray}
\left(\mathbf{E}\left[X_{\tau +1}^{2} | \mathcal{F}_{\tau}\right]\right)^p \leq \mathbf{E}\,[|X_{\tau +1}|^{2p} | \mathcal{F}_{\tau}].
\end{eqnarray}
Taking expectations on both sides of (\ref{aaa}), we deduce that
\begin{eqnarray}\label{b}
\mathbf{E}\left[\left(\sum_{i=\tau+1}^{N} \mathbf{E}\left[\big(X_{i}-\hat{X}_{i}\big)^{2} \Big| \mathcal{F}_{i-1}\right]\right)^{p}\right] &\leq&2^{p-1}s_{n}^{2p}\left\|V_{n}^{2}-1\right\|_{p}^{p}+2^{2p-1}\mathbf{E}\,[|X_{\tau +1}|^{2p}] \nonumber \\
&\leq&2^{p-1}s_{n}^{2p}\left\|V_{n}^{2}-1\right\|_{p}^{p}+2^{2p-1}\mathbf{E} \max _{1 \leq i \leq n}\left|X_{i}\right|^{2p}.
\end{eqnarray}
By an argument similar to that of (\ref{aaa}), we get
\begin{eqnarray}
\mathbf{E}\left[\max _{\tau+1 \leq i \leq N}\left|{X}_{i}-\hat{X}_{i}\right|^{2p}\right]
 &\leq& 2^{2p-1}\mathbf{E}\left[\max _{\tau+1 \leq i \leq n}\left|X_{i}\right|^{2 p}+\varepsilon^{2p}\right] \nonumber \\
 &\leq& 2^{2p-1} \mathbf{E}\left[\max _{1 \leq i \leq n}\left|X_{i}\right|^{2 p}+\varepsilon^{2p}\right].  \label{aads2ads}
\end{eqnarray}
Combining  the inequalities (\ref{aa}), (\ref{b}) and (\ref{aads2ads}) together, we have
\begin{eqnarray}
\mathbf{E}\left[|\hat{S}_{N}-S_{n}|^{2p}\right] \leq C_{p}\left(s_{n}^{2p}\left\|V_{n}^{2}-1\right\|_{p}^{p}+\mathbf{E}\max _{1 \leq i \leq  n}\left|X_{i}\right|^{2 p}+\varepsilon^{2 p}\right).
\end{eqnarray}
Now we combine this with the inequalities (\ref{1}) and (\ref{a}) together, and then let $\varepsilon \rightarrow 0$, we get
\begin{eqnarray}
 d_{w}\big(S_{n} / s_{n}\big)\leq C \frac{1}{s_{n}}\Big(\sum_{i=1}^{n} \mathbf{E}\left|X_{i}\right|^{3}\Big)^{1/3}
+C x^{1-2p}\Big(\left\|V_{n}^{2}-1\right\|_{p}^{p}+\frac{1}{s_{n}^{2p}}\mathbf{E} \max _{1 \leq i \leq  n}\left|X_{i}\right|^{2p}\Big)+2x.
\end{eqnarray}
Putting $x=\Big(\left\|V_{n}^{2}-1\right\|_{p}^{p}+\frac{1}{s_{n}^{2p}}\mathbf{E} \max _{1 \leq i \leq  n}\left|X_{i}\right|^{2p}\Big)^{1/2p},$ the final bound follows.\hfill\qed

\subsection*{Acknowledgements}
The authors would like to thank two anonymous referees and the editor for their valuable comments.
This work has been partially supported by the National Natural Science Foundation
of China (Grant nos.\,11601375 and 11971063).

\end{document}